\documentclass[a4paper,oneside,12pt]{article}
\usepackage{amsmath,amsfonts,amscd,amssymb}
\usepackage{longtable,geometry}
\usepackage[english]{babel}
\usepackage[utf8]{inputenc}
\usepackage[T1]{fontenc}
\newtheorem{theorem}{Theorem}[section]

\newtheorem{lemma}[theorem]{Lemma}

\numberwithin{equation}{section}

\newcommand{\IP}{\mathbb{P}}
\newcommand{\IE}{\mathbb{E}}
\newcommand{\B}{\mathcal{B}}
\newcommand{\Zd}{\mathbb{Z}^d}
\newcommand{\Var}{\text{Var}}
\newcommand{\Rq}{\textbf{Remark}}

\title{Asymptotics of the odometer function for the internal Diffusion Limited Aggregation model}
\author{Cyrille Lucas\footnote{Mod\'elisation al\'eatoire de Paris 10 (MODAL'X)-Universit\'e Paris Ouest Nanterre La
D\'efense, \newline 200 avenue de la République, 92001 Nanterre Cédex \newline \textit{email:} cyrille.lucas@u-paris10.fr}}
\date{November 16, 2009}

\begin{document}
\maketitle
\small\textsc{Abstract:} We present precise asymptotics of the odometer function for the internal Diffusion Limited Aggregation model. These results provide a better understanding of this function whose importance was demonstrated by Levine and Peres \cite{Peres}. We derive a different proof of a time-scale result by Lawler, Bramson and Griffeath \cite{LBG}.
 \normalsize
\begin{center} \textbf{Keywords:} asymptotic shape, divisible sandpile, Green's function, internal diffusion limited aggregation.
\end{center}

\section{Introduction}
The internal Diffusion Limited Aggregation model was first introduced by Diaconis and Fulton in \cite{DF} gives a protocol for building a random set recursively. At each step, the first vertex visited outside the cluster by a simple random walk started at the origin is added to the cluster. The resulting limit shape is the euclidian ball, as proved in 1992 by Lawler, Bramson and Griffeath in \cite{LBG}.

More recently, however, Levine and Peres \cite{Peres}, \cite{Peres2} have shown that this model is related to the rotor-router and divisible sandpile models. In the former, random walkers are replaced by eulerian walkers. In the latter, each vertex can hold 1 unit of mass, and the excess is divided equally among its neighbors when the vertex topples. Thus an initial mass at the origin becomes a stable shape after a suitable infinite series of topplings.

For each of these three models, one can define an odometer function, which will be the total number of times a walker passes through a given point (counting multiple passages of the same walker) for the internal DLA and rotor-router model, and the total mass emitted from a given point in the construction of the cluster for the divisible sandpile model. This function plays a comparable role in all three models and turns out to be instrumental in their relation. While the limiting shape of these models is known to be the Euclidian ball, the behavior of their odometer functions in the particular case where all the mass is started at the origin remains to be studied.

This paper provides a closer look at the odometer function in the case of the internal DLA model, with an almost sure convergence of the normalized functions and asymptotics of this function near the origin. These results provide in turn a new proof of the time scale of the cluster introduced in \cite{LBG}.

\section{Definitions, main results}
Let $ (S^j)_{j \in \mathbb{N}}$ a sequence of independent simple random walks on $\Zd$, and let us define the cluster $A(n)$ and stopping times $(\sigma_k)_{k \in \mathbb{N}}$ recursively in the following way:
\begin{eqnarray}
\sigma_0 & = & 0 \text{,}\nonumber\\
A(0)     & = & \{ 0 \} = \{ S^0(\sigma_0) \} \text{, and for all $j >0$,}\nonumber\\
\forall j >0, \text{ } \sigma_j & = & \inf \{ t \geq 0 : S^j(t) \not \in A(j-1) \} \text{,}\nonumber\\
                          A(j)  & = & A(j-1) \cup \{ S^j(\sigma_j) \} \text{.}\nonumber
\end{eqnarray}

Let $||.||$ denote the Euclidean norm on $\mathbb{R}^d$, $\B_r$ the Euclidean ball of radius $r$ of $\Zd$, and $\omega_d$ the volume of the unit ball of $\mathbb{R}^d$. We will consider the cluster $A(\omega_d n^d)$, which has the same volume as $\B_n$. Lawler, Bramson and Griffeath proved in \cite{LBG} that the normalized cluster $\frac{1}{n}A(\left\lfloor \omega_d n^d\right\rfloor)$ converges to the Euclidean unit ball with probability one. Lawler then improved this result, defining the inner and outer errors as follows:
\begin{eqnarray}
\delta_I(n) & = & n - \inf_{z \not \in A(\left\lfloor \omega_d n^d\right\rfloor)} ||z|| \text{,}\nonumber\\
\delta_O(n) & = & \sup_{z \in A(\left\lfloor \omega_d n^d\right\rfloor)} ||z|| -n \text{.}\nonumber
\end{eqnarray}

to get the following bounds, with probability $1$:
\begin{eqnarray}
\lim_{n \rightarrow \infty} \frac{\delta_I(n)}{n^{1/3}(\ln n)^{2}}&=&0 \text{, and} \nonumber\\
\lim_{n \rightarrow \infty} \frac{\delta_O(n)}{n^{1/3}(\ln n)^{4}}&=&0 \text{.}\label{lawler}
\end{eqnarray}

An important part in the study of this model is played by Green's functions for the simple random walk, which we introduce now, as defined in \cite{IRW}. Let $S^x$ be a simple random walk started from a point $x \in \Zd$, and $\xi_n$ be the time at which $S$ exits $\B_n$, that is :
$$\xi_n = \inf\{t \text{ : } S_t \notin \B_n\}.$$
Then we define for all $y \in \Zd$
$$G_n(x,y) = \sum_{t=0}^{\xi_n} \textbf{1}_{S^x(t)=y}.$$
This is zero if either $x$ or $y$ lies outside $\B_n$. For $d \geq 3$ we define, for all $x,y \in \Zd$,
$$G(x,y) = \sum_{t=0}^{\infty} \textbf{1}_{S^x(t)=y}.$$

Let $z \in \mathbb{R}^d$ be a non-zero point in the open unit ball. Let us write $z=(x_1, \cdots, x_d)$

Let us define $z_n := \left(\lfloor n x_1 \rfloor, \cdots  \lfloor n x_d \rfloor       \right)$ .

The odometer function $u_n(z)$ at rank $n$ measures the number of walkers passing through $z_n$ in the process of building the cluster $A(\left\lfloor \omega_d n^n\right\rfloor)$. It is defined as follows:
$$ u_n(z) := \sum_{i=0}^{\left\lfloor \omega^d n^d \right\rfloor} \sum_{t=0}^{\sigma_i} \mathbf{1}_{S^i(t) = z_n} $$

Our results are as follows:
\begin{theorem}
For all non-zero points of the open unit ball $z$, when $n$ goes to infinity,
\begin{align}
\frac{u_{n}(z)}{n^2}  \rightarrow & ||z||^2 -1 - 2\ln(||z||) && \text{almost surely if }d=2, \nonumber\\
\frac{u_{n}(z)}{n^2}  \rightarrow & ||z||^2 + \frac{2}{(d-2)||z||^{d-2}} - \frac{d}{d-2} && \text{almost surely if }d \geq 3.\label{f}
\end{align}
\label{theo}
\end{theorem}

\Rq \textbf{:} In \cite{LBG}, the authors estimate the time it takes to build a cluster of radius $n$, that is to say the total number of steps done by random walks during the construction of the cluster. The authors count these steps by estimating the number of steps for each given random walk. Our convergence result \ref{theo} allows us to take a different perspective on the problem, and count the total number of steps as the sum over all points in the cluster of the number of steps through this point.

This leads to a different proof of the result of Lawler, Bramson and Griffeath on the time scaling of the cluster, which is presented in section \ref{time}.

The functions in \ref{f} vanish when $||z||$ tends to $1$, and tend to infinity when $||z||$ tends to $0$. To understand this behavior around the origin, we can use a different scaling. Consider a sequence of non-zero points $y_n$ such that $\frac{y_n}{n}$ converges to $0$ as $n$ tends to infinity. If we now consider the following sequence of values takes by the $n$-th odometer at point $\frac{y_n}{n}$:
$$ u_n\left(\frac{y_n}{n}\right) = \sum_{i=0}^{\left\lfloor \omega^d n^d \right\rfloor} \sum_{t=0}^{\sigma_i} \mathbf{1}_{S^i(t) = y_n} $$

We get the following result:
\begin{theorem}
\label{theo2}
 Let $y_n$ be a sequence of non-zero points such that $\frac{y_n}{n}$ converges to $0$.
\begin{itemize}
\item{If $y_n$ converges to a point $a \in \Zd$, $d \geq 3$, when $n$ tends to infinity,
$$ u_{n}(\frac{y_n}{n}) \sim \text{ }\omega_d n^d G(0,a) \text{ almost surely}$$}
\item{If $||y_n||$ tends to infinity, when $n$ goes to infinity,
\begin{align}
\frac{u_{n}(\frac{y_n}{n})}{n^2}   \sim & \text{ }4 \ln\left(\frac{n}{||y_n||}\right)              && \text{almost surely, if }d=2, \nonumber\\
\frac{u_{n}(\frac{y_n}{n})}{n^2}   \sim & \text{ }\frac{2}{d-2}\left(\frac{||y_n||}{n}\right)^{2-d}  && \text{almost surely, if }d \geq 3. \nonumber
\end{align}}
\end{itemize}

\end{theorem}

\Rq \textbf{:} The question of the asymptotics of the odometer function near the boundary remains open, and is presumably linked to the difficult problem of the fluctuations of the cluster around its limiting shape.

\bigskip

Let us comment further on Theorem \ref{theo}, and give a heuristic of it when $d\geq 3$ based on the work of Levine and Peres (\cite{Peres}) (the two-dimensional case is similar if a little more technical). In the case of the divisible sandpile model, each site of $\Zd$ contains a continuous amount of mass. A site with an amount greater than $1$ can topple, that is to say keep mass $1$ and distribute the rest equally between its neighbors. With any sequence of topplings that topples each full site infinitely often, the mass approaches a limiting distribution, which does not depend on the sequence.

The odometer function $u$ for this model is defined as the total mass emitted from a given point. Since each neighbor $y$ of a given point $x$ divides mass equally between its $2d$ neighbors, the total mass received by $x$ is $\frac{1}{2d}\sum_{y \sim x} u(y)$. If we define the discrete Laplacian $\Delta$ as $\Delta f(x) = \frac{1}{2d}\sum_{y \sim x} f(y)-f(x)$, we get that
$$\Delta u(x) = \nu(x) - \sigma(x),$$
where $\sigma$ and $\nu$ are the initial and final amounts of mass at $x$, respectively.

In our case, the initial mass $\omega_d n^d$ is concentrated at the origin, and the final mass is $1$ in each fully occupied site of the cluster. Hence:
\begin{equation*}
\Delta u_n(x) =\begin{cases}
n-1& \text{if $x=0$},\\
-1& \text{if $x$ is a non-zero point in the cluster,}\\
0& \text{outside the cluster.}
\end{cases}
\end{equation*}

To solve this equation, we first introduce a function $\gamma_n$ that has discrete Laplacian $\sigma_n-1$, in our case:
$$ \gamma_n(z) = -||z||^2-\omega_d n^d G(0,z).$$

Then the odometer function $u_n$ is given by lemma 3.2 of \cite{Peres2} as the difference between the least superharmonic majorant of $\gamma_n$ and $\gamma_n$:
$$u_n = s_n - \gamma_n,$$
where $s_n = \inf\{f(x) \text{ such that $f$ is superharmonic and }f\geq \gamma_n\}$.

Recalling the notations $z=(x_1, \cdots, x_d)$ and $z_n = \left(\lfloor n x_1 \rfloor, \cdots  \lfloor n x_d \rfloor       \right)$, we get:
$$\frac{1}{n^2}\gamma_n(z_n) \rightarrow \gamma(z):=-||z||^2-\frac{2}{d-2}||z||^{2-d}$$
The function $\gamma$ is radial, and its particular shape makes it easy to determine its superharmonic majorant $s$, as one can see on Figure 1:
$$s(z)=\begin{cases}
-\frac{2}{d-2}& \text{if $||z||\leq 1$},\\
\gamma(z)& \text{if $||z|| \geq 1$.}
\end{cases}$$
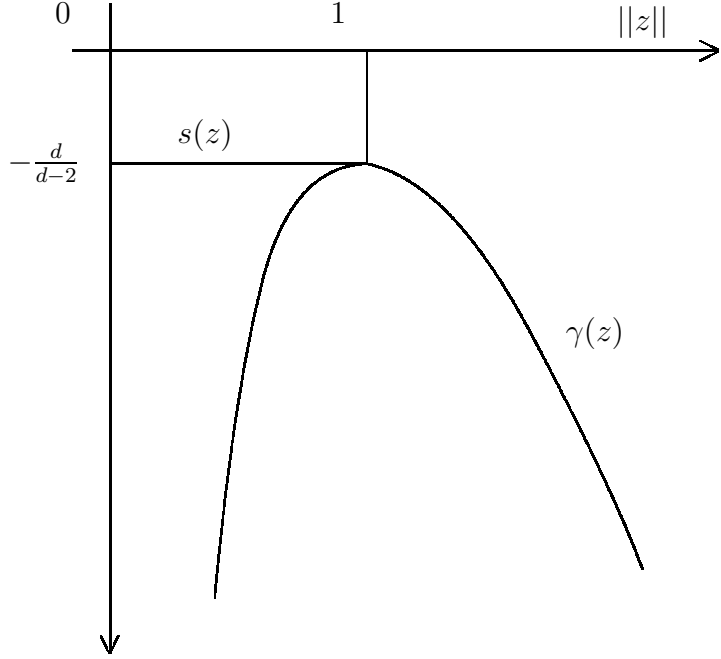
\begin{figure}
\begin{center}
\ifx\JPicScale\undefined\def\JPicScale{1}\fi
\unitlength \JPicScale mm
\begin{picture}(90,86.25)(0,0)
\linethickness{0.3mm}
\put(10,0){\line(0,1){86.25}}
\linethickness{0.3mm}
\put(5,80){\line(1,0){85}}
\linethickness{0.25mm}
\put(10,65){\line(1,0){33.75}}
\linethickness{0.05mm}
\put(43.75,65){\line(0,1){15}}
\put(3.75,85){\makebox(0,0)[cc]{$0$}}
\put(40,85){\makebox(0,0)[cc]{$1$}}
\put(1.25,65){\makebox(0,0)[cc]{$-\frac{d}{d-2}$}}
\put(22.5,68.75){\makebox(0,0)[cc]{$s(z)$}}
\put(73.75,42.5){\makebox(0,0)[cc]{$\gamma(z)$}}
\linethickness{0.2mm}
\qbezier(80,11.25)(75.74,22.81)(65.16,42.81)
\qbezier(65.16,42.81)(54.57,62.81)(43.75,65)
\qbezier(43.75,65)(33.83,64.69)(30,49.38)
\qbezier(30,49.38)(26.17,34.06)(23.75,7.5)
\linethickness{0.2mm}
\multiput(87.5,81.25)(0.25,-0.12){10}{\line(1,0){0.25}}
\linethickness{0.2mm}
\multiput(87.5,78.75)(0.25,0.12){10}{\line(1,0){0.25}}
\linethickness{0.2mm}
\multiput(10,0)(0.12,0.25){10}{\line(0,1){0.25}}
\linethickness{0.2mm}
\multiput(8.75,2.5)(0.12,-0.25){10}{\line(0,-1){0.25}}
\put(80,83.75){\makebox(0,0)[cc]{$||z||$}}
\end{picture}
\caption{The functions $s$ and $\gamma$ as radial functions of $z$.}
\end{center}
\end{figure}

Since $\gamma_n$ converges to $\gamma$, $s_n$ should converge to $s$. This would mean that $u_n$ converges to $s-\gamma$. This method is in fact made rigorous in \cite{Peres}, with the noticeable difference that in the case considered the starting mass has a bounded density, as opposed to our case where all of the mass is started at the origin. The authors then prove that the odometer function for the divisible sandpile model is the expected value of the odometer function for the internal DLA model, which gives a heuristic proof of the following result:
\begin{align}
\mathbb{E}\left( \frac{u_{n}(z)}{n^2} \right) \rightarrow & ||z||^2 -1 - 2\ln(||z||) && \text{ if }d=2, \nonumber\\
\mathbb{E}\left( \frac{u_{n}(z)}{n^2} \right) \rightarrow & ||z||^2 + \frac{2}{(d-2)||z||^{d-2}} - \frac{d}{d-2} && \text{ if }d \geq 3.\nonumber
\end{align}
This analytic method could probably be rendered rigorous, but we prefer a more probabilistic approach based on the convergence of the cluster, which gives us an almost sure convergence result.

\section{Proof of Theorem \ref{theo}}

The main idea of this proof is to use the fact that the cluster is very similar to a disc when $n$ is large enough to get two inequalities framing $u_n(z)$ between random valuations which do not depend on the shape of the cluster.

Lawler's theorem proves that the following set inequality holds for $n$ large enough with probability one:
$$ \B_{n-n^{1/3}(\ln n)^{2}} \subset A(\left\lfloor \omega_d n^d \right\rfloor) \subset \B_{n+n^{1/3}(\ln n)^{4}}, $$

Let us define:
\begin{eqnarray}
\eta_I(n) & = & n^{1/3}(\ln n)^{2} \nonumber\\
\eta_O(n) & = & n^{1/3}(\ln n)^{4} \nonumber
\end{eqnarray}

Let us define the stopping time $\xi_k^j$ as the time at which $S^j$ leaves a ball $\B_k$. As a consequence of the set inequality, we get the following inequality on stopping times for $k$ large enough, for all $j$ such that $ \left\lfloor \omega_d (k-1)^d \right\rfloor \leq j \leq \left\lfloor \omega_d k^d \right\rfloor $, with probability one:
\begin{eqnarray}
\xi_{(k-1)-\eta_I(k-1)}^j & \leq & \sigma^j \leq \xi_{k+\eta_O(k)}^j \label{st}
\end{eqnarray}

We get the following framing of $u_n(z)$:
$$
\sum_{k=0}^{n-1} \sum_{j=\left\lfloor \omega_d k^d\right\rfloor}^{\left\lfloor \omega_d (k+1)^d\right\rfloor} \sum_{t=1}^{\xi_{k-\eta_I(k)}} \mathbf{1}_{S^i(t) = z_n}  \leq  u_n(z)  \leq  \sum_{k=1}^{n} \sum_{j=\left\lfloor \omega_d k^d\right\rfloor}^{\left\lfloor \omega_d (k+1)^d\right\rfloor} \sum_{t=1}^{\xi_{k+\eta_O(k)}^j} \mathbf{1}_{S^i(t) = z_n} \label{a}
$$

This inequality holds for $n$ large enough with probability one because $||z_n|| - n||z||$ is smaller than a constant. This ensures that the inequality (\ref{st}) will start to hold as soon as $k$ becomes greater than $\frac{n||z||}{2}$ for instance. All the other terms of the sum are asymptotically zero.

We will now show that the left and right bounds on $u_n(z)$, which we will call $u_n(z)^-$ and $u_n(z)^+$, once normalized, converge almost surely to the same function. In order to do this, we take a look at the following family of random variables:
$$\chi_{k,n}^j(z) =  \sum_{t=1}^{\xi_{k}} \mathbf{1}_{S^j(t) = z_n} $$

A direct application of the Markov property shows that for $i \geq 1$, there exists $p_0$ and $p$ depending on $k$, $n$ and $z$ such that:
$$ \IP(\chi_{k,n}^j(z) =i) = p_0 p^{i-1} (1-p) $$

Since we know that $\IP(\chi_{k,n}^j(z) = 0)= \IP(\tau_{z_n} > \xi_k) $ where $\tau_{z_n}$ is the reaching time of $z_n$, and that $\IE(\chi_{k,n}^j(z)) = G_k(0,z_n)$, where $G_n$ is Green's stopped function, a simple computation shows that:
\begin{eqnarray}
p_0 & = & \IP( \xi_k < \tau_{z_n} ) \nonumber\\
p   & = & 1 - \frac{\IP( \xi_n < \tau_{z_n} )}{G_k(0,z_n)}\nonumber
\end{eqnarray}

This determines the parameter of the geometric distribution $\chi_{k,n}^j(z)$ follows.

We will need the following lemma which gives an estimate of $G_n(0,z)$:
\begin{lemma}
\label{Green1}
If $z \in \B_n$, $z \neq 0$, we have
\begin{eqnarray}
G_n(0,z) & = & \frac{2}{\omega_2} \ln \frac{n}{||z||} + O(\frac{1}{||z||}) + O(\frac{1}{n}) \text{ if $d=2$, and}\nonumber\\
G_n(0,z) & = & \frac{2}{d-2} \frac{1}{\omega_d}(||z||^{2-d} - n^{2-d}) + O(||z||^{1-d}) \text{,} \nonumber
\end{eqnarray}
where the $O$ are uniform on the unspecified variables.
\end{lemma}
This lemma is due to Lawler and can be found in section 1.6 of \cite{IRW}.

We will now prove that $\frac{1}{n^2} u_n(z)$ converges almost surely to the limit of its mean value. We first compute the mean value of $u_n(z)^-$, then we will bound its variance.

Consider the left side of (\ref{a}):
\begin{eqnarray}
u_n(z)^- & = & \sum_{k=0}^{n-1} \sum_{j=\left\lfloor \omega_d k^2\right\rfloor}^{\left\lfloor \omega_2 (k+1)^2\right\rfloor} \sum_{t=1}^{\xi_{k-\eta_I(k)}} \mathbf{1}_{S^j(t) = z_n} \nonumber\\
u_n(z)^- & = & \sum_{k=0}^{n-1} \sum_{j=\left\lfloor \omega_d k^2\right\rfloor}^{\left\lfloor \omega_2 (k+1)^2\right\rfloor} \chi_{k-\eta_I(k),n}^j(z) \nonumber\\
\IE(u_n(z)^-) & = & \sum_{k=0}^{n-1} \sum_{j=\left\lfloor \omega_d k^2\right\rfloor}^{\left\lfloor \omega_2 (k+1)^2\right\rfloor} G_{k-\eta_I(k)}(0,z_n) \nonumber\\
 & = & \sum_{k=\left\lfloor ||z||n \right\rfloor}^{n-1} (\left\lfloor \omega_d (k+1)^2\right\rfloor - \left\lfloor \omega_2 k^2\right\rfloor) G_{k-\eta_I(k)}(0,z_n) \nonumber\\
\end{eqnarray}
Let us consider separately the case $d=2$. In that case, we have:
\begin{eqnarray}
\IE(u_n(z)^-) & = & \sum_{k=\left\lfloor ||z||n \right\rfloor}^{n-1} (\left\lfloor \omega_2 (k+1)^2\right\rfloor - \left\lfloor \omega_2 k^2\right\rfloor)  \frac{2}{\omega_2} \left(\ln \frac{k-\eta_I(k)}{n||z||} + O(\frac{1}{n}) \right)\nonumber\\
\IE(u_n(z)^-) & = & \sum_{k=\left\lfloor ||z||n \right\rfloor}^{n-1}\left( 4 k \ln \frac{k}{n||z||} \right)+ O(n^{5/3}) \nonumber\\
\frac{1}{n^2}\IE(u_n(z)^-) & = &  \frac{1}{n}\sum_{k=\left\lfloor ||z||n \right\rfloor}^{n-1}\left( 4 \frac{k}{n} \ln \frac{k}{n||z||} \right)+ O(n^{-1/3}) \nonumber\\
\end{eqnarray}
We recognize a Riemann sum, thus :
\begin{eqnarray}
\frac{1}{n^2} \IE(u_n(z)^-) & \rightarrow & \int_{||z||}^{1} 4 t \ln(\frac{t}{||z||}) \text{d}t = ||z||^2 -1 - 2\ln(||z||) \nonumber
\end{eqnarray}
When $d\geq3$, we get instead the following estimate for $\IE(u_n(z)^-)$ :
$$\sum_{k=\left\lfloor ||z||n \right\rfloor}^{n-1} (\left\lfloor \omega_d (k+1)^d\right\rfloor - \left\lfloor \omega_d k^d\right\rfloor) \frac{2}{d-2} \frac{1}{\omega_d}((n||z||)^{2-d} - (k-\eta_I(k))^{2-d}) + O(n^{1-d})$$
\begin{eqnarray}
\IE(u_n(z)^-) & = & \frac{1}{n}\sum_{k=\left\lfloor ||z||n \right\rfloor}^{n-1} \left( \frac{2d}{d-2} k^{d-1} ((n||z||)^{2-d} - (k)^{2-d}) \right) + O(n^{8/3-d}) \nonumber\\
\frac{1}{n^2}\IE(u_n(z)^-) & = &  \frac{1}{n}\sum_{k=\left\lfloor ||z||n \right\rfloor}^{n-1}\left( \frac{2d}{d-2}  ((\frac{k}{n})^{d-1}||z||^{2-d} - \frac{k}{n}) \right) + O(n^{2/3-d}) \nonumber\\
\frac{1}{n^2} \IE(u_n(z)^-) & \rightarrow & \int_{||z||}^{1} \frac{2d}{d-2} ( t^{d-1}||z||^{2-d} -t ) \text{d}t = ||z||^2 + \frac{2}{(d-2)||z||^{d-2}} - \frac{d}{d-2} \nonumber
\end{eqnarray}

The same arguments show that $\frac{1}{n^2} \IE(u_n(z)^+)$ has the same limit as its counterpart.

We will now give a bound on the variance of these variables. Let us first state the following lemma:
\begin{lemma}
\label{var}
The random valuations $\chi_{k,n}^j(z)$ have the following variance, for all $j,k,n \in \mathbb{N},z\in \Zd$ provided that $||z||n \leq k$:
$$\Var (\chi_{k,n}^j(z)) = \left(1 - \frac{\IP( \xi_k < \tau_{z_n} )}{G_k(0,z_n)}\right)G_k(0,z_n)G_k(z_n,z_n)\left(2+\frac{1}{G_k(z_n,z_n)-1}\right) - G_k(0,z_n)^2 $$
\end{lemma}

The proof of this lemma is a straightforward calculation which relies only on the knowledge of the law of $\chi_{k,n}^j(z)$.

We will need bounds for $G_k(z_n,z_n)$. Provided that $||z||n \geq k$, $z_n$ is in the ball of radius $k$ centered at the origin, and it has at least one neighbor in this same ball, which proves that $G_k(z_n,z_n) \geq 1 + \frac{1}{4d^2} $.

When $d \geq 3$, we can use the simple upper bound given by $G_k(z_n,z_n) \leq G(0,0) < \infty$.

However, when $d=2$, we will use the following bound: $G_k(z_n,z_n) \leq G_{n+k}(0,0)$. This bound comes from the fact that a random walk started in $z_n$ will exit the ball of radius $k$ centered in $0$ before it exits the ball of radius $n+k$ centered in $z_n$. When $k \leq n$, we will even use $G_k(z_n,z_n) \leq G_{2n}(0,0)$.

We are now ready to bound the variance of $u_n(z)^-$.
\begin{eqnarray}
\Var \left(\frac{u_n(z)^-}{n^2}\right) & = & \frac{1}{n^4} \sum_{k=0}^{n-1} \sum_{j=\left\lfloor \omega_d k^2\right\rfloor}^{\left\lfloor \omega_2 (k+1)^2\right\rfloor} \Var(\chi_{k-\eta_I(k),n}^j(z)) \nonumber
\end{eqnarray}
$$\leq  \frac{1}{n^4} \sum_{k=0}^{n-1} \sum_{j=\left\lfloor \omega_d k^2\right\rfloor}^{\left\lfloor \omega_2 (k+1)^2\right\rfloor} G_{k-\eta_I(k)}(0,z_n)G_{k-\eta_I(k)}(z_n,z_n)\left(2+\frac{1}{G_{k-\eta_I(k)}(z_n,z_n)-1}\right)$$

When $d \geq 3$, naming $K_d$ a suited constant depending only on the dimension $d$, we get:
\begin{eqnarray}
\Var \left(\frac{u_n(z)^-}{n^2}\right) & \leq & \frac{1}{n^4} \sum_{k=0}^{n-1} K_d k^{d-1} G_k(0,z_n) \nonumber
\end{eqnarray}

This formula is the same as in the estimation of $\IE(u_n(z)^-)$ and yields:
$$ \Var \left(\frac{u_n(z)^-}{n^2}\right) \leq \frac{g(||z||,d)}{n^2} ,$$
where $g$ is a function of $||z||$ and $d$ only.

When $d = 2$, taking a suited constant $K_2$ gives us:
\begin{eqnarray}
\Var \left(\frac{u_n(z)^-}{n^2}\right) & \leq & \frac{1}{n^4} \sum_{k=0}^{n-1} K_2 k G_k(0,z_n) \ln(n)^2, \nonumber
\end{eqnarray}
which yields:
$$ \Var \left(\frac{u_n(z)^-}{n^2}\right) \leq \frac{\ln(n)^2h(||z||)}{n^2} ,$$
where $h$ is a function of $||z||$.

In both cases, the sum $\sum \Var \left(\frac{u_n(z)^-}{n^2}\right)$ is finite, which means $\frac{u_n(z)^-}{n^2}$ converges almost surely to $\lim_{n \rightarrow \infty} \frac{\IE(u_n(z)^-)}{n^2}$. The same set of arguments can be used to prove that $\frac{u_n(z)^+}{n^2}$ converges almost surely to $\lim_{n \rightarrow \infty} \frac{\IE(u_n(z)^+)}{n^2}$.

Since $\lim_{n \rightarrow \infty} \frac{\IE(u_n(z)^-)}{n^2} = \lim_{n \rightarrow \infty} \frac{\IE(u_n(z)^+)}{n^2}$ and $ u_n^-(z) \leq u_n(z) \leq u_n^+(z)$ for all $n$ large enough, we have proved our result.

\section{Proof of Theorem \ref{theo2}}

The proof of the first assertion relies only on the fact that the following framing holds for $n$ large enough:
$$
\sum_{k=\left\lfloor \ln n\right\rfloor}^{n-1} \sum_{j=\left\lfloor \omega_d k^d\right\rfloor}^{\left\lfloor \omega_d (k+1)^d\right\rfloor} \sum_{t=1}^{\xi_{k-\eta_I(k)}} \mathbf{1}_{S^i(t) = y_n}
\leq  \sum_{i=\left\lfloor \omega_d (\ln n)^d \right\rfloor}^{\left\lfloor \omega_d n^d \right\rfloor} \sum_{t=0}^{\sigma_i} \mathbf{1}_{S^i(t) = y_n}  $$
$$\leq \sum_{k=1+\left\lfloor \ln n\right\rfloor}^{n} \sum_{j=\left\lfloor \omega_d k^d\right\rfloor}^{\left\lfloor \omega_d (k+1)^d\right\rfloor} \sum_{t=1}^{\xi_{k+\eta_O(k)}^j} \mathbf{1}_{S^i(t) = y_n} \nonumber
$$

This is true because the equation ($\ref{st}$) holds as soon as $k$ is large enough, which happens eventually since $k \geq \ln n$.

The difference to the relevant quantity can be bounded using this inequality which holds for $n$ large enough:
\begin{eqnarray}
\sum_{i=0}^{\left\lfloor \omega^d (\ln n)^d \right\rfloor} \sum_{t=0}^{\sigma_i} \mathbf{1}_{S^i(t) = y_n}
& \leq & \sum_{i=0}^{\left\lfloor \omega_d (\ln n)^d\right\rfloor}  \sum_{t=1}^{\xi_{\ln n-\eta_I(\ln n)}} \mathbf{1}_{S^i(t) = y_n} \nonumber
\end{eqnarray}

It is a consequence of the calculations in the proof of Theorem \ref{theo} that the right side tends to zero almost surely once divided by $n^2$.

Just like in the proof of Theorem \ref{theo}, we estimate the expected value of our lower bound:
\begin{eqnarray}
\mathbb{E}\left( \sum_{k=\left\lfloor \ln n\right\rfloor}^{n-1} \sum_{j=\left\lfloor \omega_d k^d\right\rfloor}^{\left\lfloor \omega_d (k+1)^d\right\rfloor} \sum_{t=1}^{\xi_{k-\eta_I(k)}} \mathbf{1}_{S^i(t) = y_n} \right)
& = & \sum_{k=0}^n d \omega_d k^{d-1} G(||y_n||) + O(n^2) \nonumber\\
& = & \omega_d n^d G(a) + O(n^2) \nonumber
\end{eqnarray}

The calculation for the upper bound yields the same result.

To finish the proof, we just need to bound the variance of these two random variables, using lemma \ref{var}:
$$\Var\left(u_n\left(\frac{y_n}{n}\right)^-\right)= \sum_{k=0}^{n-1} \sum_{j=\left\lfloor \omega_d k^2\right\rfloor}^{\left\lfloor \omega_d (k+1)^2\right\rfloor} \Var(\chi_{k-\eta_I(k),n}^j(\frac{y_n}{n}))$$
\begin{eqnarray}
  &  \leq & \sum_{k=0}^{n-1} \sum_{j=\left\lfloor \omega_d k^2\right\rfloor}^{\left\lfloor \omega_d (k+1)^2\right\rfloor} G_{k-\eta_I(k)}(0,y_n)G_{k-\eta_I(k)}(y_n,y_n)\left(2+\frac{1}{G_{k-\eta_I(k)}(y_n,y_n)-1}\right) \nonumber
\end{eqnarray}

When $d \geq 3$, naming $K_d$ a suited constant depending only on the dimension $d$, we get:
\begin{eqnarray}
\Var \left( u_n\left(\frac{y_n}{n}\right) \right)^- & \leq & \sum_{k=0}^{n-1} K_d k^{d-1} G(||y_n||) \nonumber
\end{eqnarray}

Hence we get:
\begin{eqnarray}
\Var \left(\frac{u_n\left(\frac{y_n}{n}\right)^-}{n^d}\right) & \leq & \frac{K'_d G(a)}{n^d} \nonumber
\end{eqnarray}

When $d = 2$, taking a suited constant $C_2$ gives us:
\begin{eqnarray}
\Var \left( u_n\left(\frac{y_n}{n}\right) \right) & \leq & \sum_{k=0}^{n-1} K_2 k G(||y_n||) \ln(n)^2, \nonumber
\end{eqnarray}
which yields:
$$ \Var \left(\frac{u_n(z)^-}{n^2}\right) \leq \frac{\ln(n)^2 K'_2 G(a)}{n^2} ,$$
where $C'_2$ is a suited constant.

In both cases, the sum $\sum \Var \left(\frac{u_n\left(\frac{y_n}{n}\right)^-}{n^d}\right)$ is finite. Since the same calculation can be applied to $u_n^+$, we can use the same argument as in the proof of Theorem \ref{theo} to say that almost surely, as $n$ tends to infinity,
$$ \frac{u_{n}(\frac{y_n}{n})}{n^d} \rightarrow \omega_d G(a). $$

The proof of the second assumption is similar to that of the first one and relies on the following estimate: for $d=2$,
\begin{eqnarray}
\mathbb{E}\left(\sum_{k=0}^{n-1} \sum_{j=\left\lfloor \omega_2 k^2\right\rfloor}^{\left\lfloor \omega_2 (k+1)^2\right\rfloor} \sum_{t=1}^{\xi_{k-\eta_I(k)}} \mathbf{1}_{S^i(t) = y_n} \right)
& = & \sum_{k=0}^n 4 k \ln\left(\frac{||y_n||^{d-2}}{k-\eta_I(k)}\right)+ O\left( \frac{1}{||y_n||} + \frac{1}{n} \right) \nonumber\\
& = & 4n^2 \ln\left(\frac{n}{||y_n||}\right) + O\left(\frac{1}{||y_n||} + \frac{1}{n} \right) \nonumber
\end{eqnarray}

And for $d\geq 3$,
\begin{gather}
\begin{split}
\mathbb{E}\left(\sum_{k=0}^{n-1} \sum_{j=\left\lfloor \omega_d k^d\right\rfloor}^{\left\lfloor \omega_d (k+1)^d\right\rfloor} \sum_{t=1}^{\xi_{k-\eta_I(k)}} \mathbf{1}_{S^i(t) = y_n} \right)
& =   \sum_{k=0}^n \frac{2d}{d-2} k^{d-1} \left( \frac{1}{||y_n||^{d-2}} - \frac{1}{(k-\eta_I(k))^{d-2}} \right)\nonumber\\
& + O\left( \frac{n^{d-1}}{||y_n||^{d-1}}\right) \nonumber
\end{split}
\end{gather}
$$ =  \frac{2n^2}{d-2} \left( \left( \frac{n}{||y_n||} \right)^{d-2} -\frac{1}{2} \right) + O\left(\frac{n^{d-1}}{||y_n||^{d-1}} + \frac{1}{n}\right) \nonumber$$

The computations of the variances yield, for $d=2$,
$$\Var\left( u_n\left(\frac{y_n}{n} \right) \right) \leq J_2 n^2 \left( \ln n \right)^2 \ln\left(\frac{n}{||y_n||}\right),$$

and for $d\geq 3$,
$$\Var\left( u_n\left(\frac{y_n}{n} \right) \right) \leq J_d n^2 \left(\frac{||y_n||}{n}\right)^{2-d},$$
where $J_2$ and $J_d$ are suitable constants depending only on $d$. It follows that
\begin{eqnarray}
\Var\left(\frac{\left(u_n\left(\frac{y_n}{n} \right)\right)}{n^2 \ln \left( \frac{n}{||y_n||} \right)}\right) &\leq& \frac{J_2 \ln^2(n)}{n^2 \ln\left(\frac{n}{||y_n||}\right)} \nonumber\\
\Var\left(\frac{u_n\left(\frac{y_n}{n} \right)}{n^2 \ln \left( \frac{n}{||y_n||} \right)}\right) &\leq&  \frac{J_d}{n^2} \left(\frac{||y_n||}{n}\right)^{d-2} \nonumber
\end{eqnarray}

These quantities have finite sums, which proves our result.

\section{Back to a result of Lawler, Bramson and Griffeath}
\label{time}
The result in Theorem \ref{theo} can be used to derive a different proof of an existing time-scale result by Lawler, Bramson and Griffeath.

In \cite{LBG}, the authors are interested in the time it takes to build a cluster of radius $n$, that is to say the total number of steps done by random walks during the construction of the cluster. They define $\tilde A(t)$ as the internal DLA cluster on the time scale of individual random walks. Namely, if
$$ \chi(t) = \max\{k \text{ such that } \sigma_1 + \cdots + \sigma_k \leq t \}$$

then $\tilde A (t) = A(\chi(t))$.

\begin{theorem}[Lawler, Bramson, Griffeath, 1992]
Let $\tilde A(t)$ be the internal DLA cluster on the time scale of individual random walks, then for all $\epsilon > 0$, almost surely,
$$ \B_{n(1-\epsilon)} \subset \tilde A(t_n) \subset \B_{n(1+\epsilon)}, $$

where $$t_n = n^{d+2}\frac{d\omega_d}{d+2}.$$
\end{theorem}

While the original proof of this theorem studies the time spent by a given random walk inside the cluster, we have studied the time spent in a given point by all the random walks. This leads naturally to a different approach of the problem, in which we will sum the odometer function over all the points in the cluster.

Let us call $t'_n$ the time taken to build the cluster $A(n)$. Then
$$ t'_n = \sum_{z : nz \in \mathbb{Z}^d \cap A(n)} u_n(z) $$

Then for all $n$,
$$ \sum_{z , \text{ $nz$ }\in \B_{n-\eta_{I}(n)}} u_n(z) \leq t'_n \leq \sum_{z , \text{ $nz$ }\in \B_{n+\eta_O(n)}} u_n(z) $$

We can compute the following asymptotics, as $n$ tends to infinity:
\begin{eqnarray}
\sum_{z , \text{ $nz$ }\in \B_{n-\eta_{I}(n)}} u_n(z) & = & \sum_{k=1}^n d \omega_d k^{d-1} u_n\left(\frac{k}{n}\right) + o(n^{d+2})\nonumber\\
& = & d \text{ } \omega_d n^2 \sum_{k=1}^n k^{d-1} f\left(\frac{k}{n}\right) + o(n^{d+2})\text{, almost surely.} \nonumber
\end{eqnarray}
where $f$ is defined by (\ref{f}). Hence:
\begin{eqnarray}
\sum_{z : nz \in \B_{n-\eta_{I}(n)}} u_n(z) & = & d \text{ }\omega_d n^{d+2} \int_0^1 x^{d-1} f(x) dx + o(n^{d+2}) \nonumber\\
                                                & = & n^{d+2}\left(\frac{d\text{ }\omega_d}{d+2}\right) +o(n^{d+2}) \nonumber
\end{eqnarray}

The same asymptotics hold for $\sum_{z , \text{ $nz$ }\in \B_{n+\eta_{O}(n)}} u_n(z)$, so that we have, almost surely,
$$\lim_{n \rightarrow \infty} \frac{t'_n}{n^{d+2}} = \frac{d\text{ }\omega_d}{d+2}.$$

\end{document}